\documentclass[preprint,12pt]{elsarticle}
\usepackage{amssymb}
\usepackage{enumerate}
\usepackage{amsmath,amssymb}
\newtheorem{thm}{Theorem}
\newtheorem{df}{Definition}

\newtheorem{prop}[thm]{Proposition}
\newtheorem{cor}[thm]{Corollary}
\newdefinition{rmk}{Remark}
\newproof{pf}{\textbf{Proof}}

\begin{document}

\begin{frontmatter}
\title{Some properties of operator-valued frames}
\author{L. G\u avru\c ta, P. G\u avru\c ta}
\ead {gavruta\_laura@yahoo.com, pgavruta@gmail.com}
\address{"Politehnica" University of Timi\c soara, Department of Mathematics, Pia\c{t}a Victoriei no.2, 300006 Timi\c{s}oara, Romania}

\begin{abstract}
Operator-valued frames (or $g$-frames) are generalizations of frames and fusion frames and have been used in packets encoding, quantum computing, theory of coherent states and more. In this paper, we give a new formula for operator-valued frames for finite dimensional Hilbert spaces. As an application, we derive in a simple manner a recent result of A. Najati concerning the approximation of $g$-frames by Parseval ones. We obtain also some results concerning the best approximation of operator-valued frames by its alternate duals, with optimal estimates.
\end{abstract}


\end{frontmatter}

\section{Introduction}

Frames in Hilbert spaces were introduced by Duffin and Schaeffer \cite{Duffin} in 1952, in the context of nonharmonic Fourier series. After a couple of years, in 1986, frames were brought to life by Daubechies, Grossman and Meyer \cite{Daubechies}.
 Frames have nice properties which makes them useful tools in  signal processing, image processing, coding theory, sampling theory and more


In the following we denote by $\mathcal{H}$ a separable Hilbert space and by
$\mathcal{L}(\mathcal{H})$ the space of all linear bounded operators on $\mathcal{H}.$
 \begin{df} A  family of elements $\{f_n\}_{n=1}^{\infty}\subset\mathcal{H}$ is called a frame for $\mathcal{H}$ if there exists
 constants $A,B>0$ such that
 $$A\|x\|^2\leq \sum_{n=1}^{\infty}|\langle x,f_n\rangle|^2\leq B\|x\|^2,\quad x\in\mathcal{H}.$$
 The constants $A,B$ are called frame bounds.
\end{df}

We say that a frame is tight if $A=B$, a Parseval frame if $A=B=1$
and an exact frame if it ceases to be a frame when any one of its elements is removed.

The exact frames are in fact Riesz bases. If just the last inequality in the above definition holds, we say that  $\{f_n\}_{n=1}^{\infty}$ is a Bessel sequence.

The operator $$T:l^2\rightarrow \mathcal{H},\hspace{2mm}T\{c_n\}_{n=1}^{\infty}:=\sum_{n=1}^{\infty} c_nf_n$$ is called synthesis operator(or pre-frame operator). The adjoint operator is given by $$\Theta=T^{*}:\mathcal{H}\rightarrow l^2,\hspace{2mm}\Theta x=\{\langle x,f_n\rangle\}_{n=1}^{\infty}.$$ and is called the analysis operator. By composing $T$ with its adjoint $T^{*}$ we obtain the frame operator
$$S:\mathcal{H}\rightarrow\mathcal{H},\hspace{2mm}Sx=TT^{*}x=\sum_{n=1}^{\infty}\langle x,f_n\rangle f_n.$$

The next theorem is one of the most important results about frames.
\begin{thm} Let $\{f_n\}_{n=1}^{\infty}\subset \mathcal{H}$ be a frame for $\mathcal{H}$ with frame operator $S$. Then
\begin{enumerate}[($i$)]
\item $S$ is invertible and self-adjoint;
\item every $x\in \mathcal{H}$ can be represented as
\begin{equation}\label{descompunerea_reperului}x=\sum_{n=1}^{\infty}\langle x,f_n\rangle  S^{-1}f_n=\sum_{n=1}^{\infty}\langle x,S^{-1}f_n\rangle f_n.
\end{equation}
\end{enumerate}
\end{thm}
The relation (\ref{descompunerea_reperului}) is called the reconstruction formula. We call $\{\langle x,S^{-1}f_n\rangle\}_{n=1}^{\infty}$ the frame coefficients.

	 The frame $\{S^{-1}f_n\}_{n=1}^{\infty}$ is called the canonical dual frame of  $\{f_n\}_{n=1}^{\infty}.$
	 A sequence $\{g_n\}_{n=1}^{\infty}$ for $\mathcal{H}$ is called an alternate dual for $\{f_n\}_{n=1}^{\infty}$ if it satisfies the following equality $$x=\sum_{n=1}^{\infty}\langle x,g_n\rangle  f_n,\quad \forall x\in\mathcal{H}.$$

A generalization of frames, which allows to reconstruct elements from the range of a linear and bounded operator in a Hilbert space, was obtain by L. G\u avru\c ta \cite{LGavruta0}.

In 2006, W. Sun \cite{Sun} introduced the concept of $g$-frame. $g$-frames are generalized frames which include ordinary frames, bounded
invertible linear operators, fusion frames, as well as many recent generalizations of frames. See also the paper of V. Kaftal, D. Larson and S. Zhang \cite{Kaftal}. For the general theory of fusion frames see the papers of P.G. Casazza et al. \cite{Casazza5} and P. G\u avru\c ta \cite{PGavruta}.

For the connection between the theory of g-frames and quantum theory as in \cite{Choi} and \cite{Kraus}, see the papers \cite{Bodmann}, \cite{Han}.

 In the following we consider $\mathcal{H}$ and $\mathcal{K}$ be two Hilbert spaces.
    We denote by $\mathcal{L}(\mathcal{H},\mathcal{K})$ the space of all linear bounded operators from $\mathcal{H}$ into $\mathcal{K}$.
 By $\mathbb{I}$ we denote a finite or a countable set.

\begin{df}We say that a sequence $\{\Lambda_i\in\mathcal{L}(\mathcal{H},\mathcal{K}):i\in\mathbb{I}\}$ is a generalized frame or a $g$-frame for $\mathcal{H}$ if there exists two positive constants $A$ and $B$ such that
$$A\|x\|^2\leq\sum_{i\in\mathbb{I}}\|\Lambda_ix\|^2\leq B\|x\|^2,\quad\forall x\in\mathcal{H}.$$
\end{df}

We call $A$ and $B$ frame bounds.
We say that $\{\Lambda_i: i\in\mathbb{I}\}$ is a \textit{$g$-tight frame} if $A=B$
and a \textit{$g$-Parseval frame} if $A=B=1.$

A frame is equivalent to a g-frame whenever $\mathcal{K}=\mathbb{C}.$

The $g$-frame operator $S$ is defined as follows $$Sx=\sum_{i\in\mathbb{I}}\Lambda_i^*\Lambda_ix,$$ where $\Lambda_i^*$ is the adjoint operator of $\Lambda_i$.

W. Sun proved in the paper \cite{Sun} that $S$ is well-defined, bounded and self-adjoint operator. Then the following reconstruction formula takes place for all $x\in\mathcal{H}$
$$x=SS^{-1}x=S^{-1}Sx=\sum_{i\in\mathbb{I}}\Lambda_i^*\Lambda_iS^{-1}x=\sum_{i\in\mathbb{I}}S^{-1}\Lambda_i^*\Lambda_ix$$
 We call $\{\Lambda_iS^{-1}\}$ the \textit{canonical dual $g$-frame} of $\{\Lambda_i\}.$
  A $g$-frame $\{\Gamma_i\}$ is called an \textit{alternate dual $g$-frame} of $\{\Lambda_i\}$ if it satisfies $$x=\sum_{i\in\mathbb{I}}\Lambda_i^*\Gamma_ix,\quad\forall x\in\mathcal{H}.$$

\section{Results}

Before the main results, we give some preliminary results (we refer to Propositions \ref{prop1}, \ref{prop2}, \ref{prop3}, \ref{propref}).
\begin{prop}\label{prop1}
Let $\{\Lambda_i\}\in\mathcal{L}(\mathcal{H}, \mathcal{K})$ be a $g$-frame. Then $\mathcal{H}$ is finite-dimensional iff $$\sum_i\|\Lambda_i\|^2_F<\infty,$$
where by $\|\cdot\|_F$ we denote the Frobenius norm (or the Hilbert-Schmidt norm).
\end{prop}
\begin{pf}Let $\{e_j\}_{j\in\mathbb{J}}$ be an orthonormal basis for $\mathcal{H}$. Then $$\sum_i\|\Lambda_i\|_F^2=\sum_i\sum_k\|\Lambda_ie_k\|^2=\sum_k\sum_i\|\Lambda_ie_k\|^2$$
If dim $\mathcal{H}=n=$ card $\mathbb{J}$ we have $\displaystyle\sum_i\|\Lambda_i\|_2^2\leq\sum_{k\in\mathbb{J}}B=B$ card
$\mathbb{J}<\infty.$\\
If $\displaystyle\sum_k\|\Lambda_i\|_F^2<\infty$ we have $A$ card $\mathbb{J}\leq\displaystyle\sum_k\sum_i\|\Lambda_ie_k\|^2<\infty.$
\end{pf}
\begin{prop}\label{prop2}If $\{\Lambda_i\}$ and $\{\Gamma_i\}$ are Parseval $g$-frames we have the following equality $$\sum_i\|L\Lambda^*_i\|^2_F=\sum_i\|L\Gamma^*_i\|^2_F,$$
for $L\in\mathcal{L}(\mathcal{H},\mathcal{K}).$
\end{prop}
\begin{pf}\begin{align*}
\sum_i\|L\Lambda^*_i\|_F^2&=\sum_i\|\Lambda_iL^*\|_F^2\\
                        &=\sum_i\sum_k\|\Lambda_iL^*e_k\|^2\\
                        &=\sum_k\sum_i\|\Lambda_i(L^*e_k)\|^2\\
                        &=\sum_k\|L^*e_k\|^2
\end{align*}
We use the fact that $\|T\|_F=\|T^*\|_F.$
\end{pf}
\begin{prop}\label{prop3} Let $\{\Gamma_i\}$ be a Parseval $g$-frame and $\mathcal{H}$ a finite n-dimensional Hilbert space. Then $$\sum_{i\in\mathbb{I}}\|\Gamma_i\|_F^2=n.$$
\end{prop}
\begin{pf}$\displaystyle \sum_i\|\Gamma_i\|_F^2=\sum_i\sum_k\|\Gamma_ie_k\|^2=\sum_k\sum_i\|\Gamma_ie_k\|^2=\sum_k 1=n.$
\end{pf}

More general, we have the following result.
\begin{prop}\label{propref}
Let  $\{\Lambda_i\}$ be a g-frame  for a finite $n$-dimensional Hilbert space $\mathcal{H}$ with g-frame operator $S$.
Then for all real number $a$, we have
$$\sum_i\|\Lambda_i S^a\|_F^2=\rm{Tr}(S^{2a+1}).$$
\end{prop}
\begin{pf}
Let $\{e_k\}_{k=1}^n$ be an orthonormal basis for $\mathcal{H}$ and $a\in\mathbb{R}$. Then
\begin{align*}
\sum_i\|\Lambda_iS^a\|_F^2&=\sum_i\|S^a\Lambda_i^*\|_F^2=\sum_i\sum_j\langle S^a\Lambda_i^*e_j,
S^a \Lambda_i^*e_j\rangle\\
&=\sum_i\sum_j\langle e_j,\Lambda_iS^{2a} \Lambda_i^*e_j\rangle=\sum_i\rm{Tr}(\Lambda_iS^{2a} \Lambda_i^* )\\
&=\sum_i\rm{Tr}(S^{2a} \Lambda_i^*\Lambda_i )=\rm{Tr}(S^{2a} \sum_i\Lambda_i^*\Lambda_i )= \rm{Tr}(S^{2a+1}).
\end{align*}
\end{pf}
The identity given in the next theorem was obtain for the particular case of ordinary frames in the paper \cite{LGavruta} and for the case of continuous frames in the paper \cite{LGavruta1}.
\begin{thm}\label{thm5} Let $\{\Lambda_i\}$ be a $g$-frame with $g$-frame operator $S$, $\{\Gamma_i\}$ be a Parseval $g$-frame and $\mathcal{H}$ be a finite-dimensional Hilbert space. Then we have the following estimation: $$\sum_i\|\Lambda_i-\Gamma_i\|^2_F=\sum_i\|\Lambda_i-\Lambda_iS^{-\frac{1}{2}}\|^2_F+\sum_i\|\Gamma_iS^{\frac{1}{4}}-\Lambda_iS^{-\frac{1}{4}}\|^2_F.$$
\end{thm}
\begin{pf}We have
\begin{align*}\sum_i\|\Lambda_i-\Gamma_i\|_F^2&-\sum_i\|\Lambda_i-\Lambda_iS^{-1/2}\|_F^2\\
&=\sum_i\|\Lambda_i^*-\Gamma_i^*\|^2_F-\sum_i\|\Lambda_i^*-S^{-1/2}\Lambda_i^*\|_F^2\\
&=-2Re\sum_i\langle\Lambda_i^*,\Gamma_i^*\rangle+2Re\sum_i\langle \Lambda_i^*, S^{-1/2}\Lambda_i^*\rangle\end{align*}
On the other hand, \begin{align*}
\sum_i\|\Gamma_iS^{1/4}&-\Lambda_iS^{-1/4}\|_F^2=\sum_i\|S^{1/4}\Gamma_i^*-S^{-1/4}\Lambda_i^*\|_F^2\\
                                                &=\sum_i\|S^{1/4}\Gamma_i^*\|_F^2+\sum_i\|S^{-1/4}\Lambda_i^*\|_F^2-2
                                                Re\sum_i\langle S^{1/4}\Gamma_i^*, S^{-1/4}\Lambda_i^*\rangle\\
                                                &=\sum_i\|S^{1/4}(S^{-1/2}\Lambda_i^*)\|_F^2+\sum_i\|S^{-1/4}\Lambda_i^*\|_F^2-2Re\sum_i\langle \Gamma_i^*,\Lambda_i^*\rangle\\
                                                &=2Re\sum_i\langle\Lambda_i^*,S^{-1/2}\Lambda_i^*\rangle-2Re\sum_i\langle\Gamma_i^*,\Lambda_i^*\rangle
\end{align*}
where we use Proposition \ref{prop3}.
\end{pf}
As a corollary, we have immediately the following result of A. Najati \cite{Najati}.
\begin{cor}Let $\{\Lambda_i\}$ be a $g$-frame, with $g$-frame operator $S$ and $\mathcal{H}$ is a finite-dimensional Hilbert space. For all Parseval $g$-frames $\{\Gamma_i\}$ the following inequality
$$\sum_i\|\Lambda_i-\Gamma_i\|^2_F\geq\sum_i\|\Lambda_i-\Lambda_iS^{-\frac{1}{2}}\|^2_F$$ takes place, and we have equality iff $\Gamma_i=\Lambda_i S^{-\frac{1}{2}}.$
\end{cor}

In the following we consider $0\leq\varepsilon<1.$ The next definition extends the definition for ordinary frames of P.G. Casazza \cite{Casazza2}.
\begin{df}
We say that $\{\Lambda_i\}$ is $\varepsilon$ nearly Parseval $g$-frame if for all $x\in\mathcal{H}$,
\begin{equation}\label{defepsilon}(1-\varepsilon)\|x\|^2\leq\sum_i\|\Lambda_ix\|^2\leq(1+\varepsilon)\|x\|^2.\end{equation}
\end{df}
\begin{thm} If $\{\Lambda_i\}$ is an $\varepsilon$ nearly Parseval $g$-frame and $\mathcal{H}$ is a finite n-dimensional Hilbert space, then we have
$$\sum_i\|\Lambda_i-\Lambda_i S^{-\frac{1}{2}}\|_F^2\leq n(1-\sqrt{1-\varepsilon})^2.$$
Moreover, the estimation is optimal.
\end{thm}
\begin{pf} Let be $S$, with eigenvalues $\{\lambda_k\}_{k=1}^n$ and a correspondent orthonormal set of eigenvectors $\{e_k\}_{k=1}^n.$ It follows that
\begin{align*}
\sum_i\|\Lambda_i-\Lambda_iS^{-1/2}\|_F^2&=\sum_i\sum_k\|\Lambda_ie_k-\Lambda_iS^{-1/2}e_k\|^2\\
                                         &=\sum_k\sum_i\|\Lambda_ie_k-\frac{1}{\sqrt{\lambda_k}}\Lambda_ie_k\|^2\\
                                         &=\sum_k\bigg(1-\frac{1}{\sqrt{\lambda_k}}\bigg)^2\sum_i\|\Lambda_ie_k\|^2\\
                                         &=\sum_k\bigg(1-\frac{1}{\sqrt{\lambda_k}}\bigg)^2\langle Se_k,e_k\rangle\\
                                         &=\sum_k\bigg(1-\frac{1}{\sqrt{\lambda_k}}\bigg)^2\lambda_k\langle e_k,e_k\rangle=\sum_k(\sqrt{\lambda_k}-1)^2
                                         \end{align*}
If in (\ref{defepsilon}), we put $x=e_k$, we obtain $$1-\varepsilon\leq\langle Se_k,e_k\rangle\leq 1+\varepsilon$$
and from here $1-\varepsilon\leq \lambda_k\leq 1+\varepsilon$, which implies $$-1+\sqrt{1-\varepsilon}\leq\sqrt{\lambda_k}-1\leq \sqrt{1+\varepsilon}-1\leq 1-\sqrt{1-\varepsilon},$$
since$\sqrt{1-\varepsilon}+\sqrt{1+\varepsilon}\leq 2$ and thus $$\sum_i\|\Lambda_i-\Lambda_iS^{-1/2}\|_F^2\leq n(1-\sqrt{1-\varepsilon})^2.$$

We prove now that the estimation is optimal. Indeed, we take $\{\Lambda_i\}=\{f_j^{\varepsilon}\}_{j=1}^n$, where $f_j^{\varepsilon}=\sqrt{1-\varepsilon}u_k,\quad k=1,2,\ldots,n$ and $\{u_k\}_{k=1}^n$ an orthonormal basis of $\mathcal{H}.$ For details see \cite{LGavruta1}
\end{pf}
In the following we give analogous results for the case when a canonical Parseval $g$-frame is replaced by canonical dual frame. We have before the following identity.
\begin{prop}\label{prop1'} Let $\{\Lambda_i\}$ be a $g$-frame, with $g$-frame operator $S$ and $\{\Gamma_i\}$ be an alternate dual $g$-frame of $\{\Lambda_i\}$. Then we have the following estimation:
\begin{equation}\label{1'}
\sum_i\|\Lambda_ix-\Gamma_ix\|^2=\sum_i\|\Lambda_ix-\Lambda_iS^{-1}x\|^2+\sum_i\|\Lambda_iS^{-1}x-\Gamma_ix\|^2
\end{equation}
\end{prop}
\begin{pf} We have \begin{equation}\label{2'}
\sum_i\|\Lambda_ix-\Gamma_ix\|^2=\sum_i\|\Lambda_ix\|^2+\sum_i\|\Gamma_ix\|^2-2\langle x,x\rangle
\end{equation}
because $$\sum_i\langle\Lambda_ix,\Gamma_ix\rangle=\sum_i\langle\Gamma_i^*\Lambda_i x,x\rangle=\langle\sum_i\Gamma_i^*\Lambda_i x,x\rangle=\langle x,x\rangle$$
From equation (\ref{2'}), by taking $\Gamma_i=\Lambda_i S^{-1}$, we have
$$\sum_i\|\Lambda_ix-\Lambda_iS^{-1}x\|^2=\sum_i\|\Lambda_ix\|^2+\sum_i\|\Lambda_iS^{-1}x\|^2-2\langle x,x\rangle.$$
Using the fact that $S$ is $g$-frame operator, i.e $Sx=\displaystyle \sum_i\Lambda_i^*\Lambda_ix$ witch implies $\displaystyle\langle Sx,x\rangle=\sum_i\|\Lambda_ix\|^2$. By putting $x\mapsto S^{-1}x$ we get
\begin{equation}\label{*}
\langle x,S^{-1}x\rangle=\sum_i\|\Lambda_iS^{-1}x\|^2
\end{equation}
Then \begin{equation}\label{3'}
\sum_i\|\Lambda_ix-\Lambda_iS^{-1}x\|^2=\sum_i\|\Lambda_ix\|^2+\langle x, S^{-1}x\rangle-2\langle x,x\rangle
\end{equation}
We also have \begin{align*}
\sum_i\|\Gamma_ix-\Lambda_iS^{-1}x\|^2&=\sum_i\|\Gamma_ix\|^2+\sum_i\|\Lambda_iS^{-1}x\|^2-2\sum_i\langle\Gamma_ix,\Lambda_iS^{-1}x\rangle\\
&=\sum_i\|\Gamma_ix\|^2+\langle x, S^{-1}x\rangle-2\sum_i\langle\Lambda_i^*\Gamma_ix, S^{-1}x\rangle
\end{align*}
So \begin{equation}\label{4'}
\sum_i\|\Gamma_ix-\Lambda_iS^{-1}x\|^2=\sum_i\|\Gamma_ix\|^2-\langle x, S^{-1}x\rangle
\end{equation}
If we add relations (\ref{3'}), (\ref{4'}) and using (\ref{2'}) we get \begin{align*}\sum_i\|\Lambda_ix-\Lambda_iS^{-1}x\|^2+\sum_i\|\Gamma_ix-\Lambda_iS^{-1}x\|^2&=
\sum_i\|\Lambda_ix\|^2+\sum_i\|\Gamma_ix\|^2-2\langle x,x\rangle\\
&=\sum_i\|\Lambda_ix-\Gamma_ix\|^2
\end{align*}
So we obtained (\ref{1'}).
\end{pf}
\begin{cor} Let $\{\Lambda_i\}$ be a $g$-frame for the Hilbert space $\mathcal{H}$, with frame
operator $S$. For all alternate dual g-frames $\{\Gamma_i\}$ of $\{\Lambda_i\}$, the inequality
$$\sum_i\|\Lambda_ix-\Gamma_ix\|^2\geq\sum_i\|\Lambda_ix-\Lambda_iS^{-1}x\|^2$$
takes place and we have equality if and only if $\Gamma_i=\Lambda_iS^{-1}.$
\end{cor}
The following result is an analog of Theorem \ref{thm5} for alternate duals.
\begin{thm} Let $\{\Lambda_i\}$ be a $g$-frame for $\mathcal{H}$ finite dimensional Hilbert space, with frame
operator $S$ and $\{\Gamma_i\}$ be an alternate dual of $\{\Lambda_i\}$. Then we have the following estimation:
$$\sum_i\|\Lambda_i-\Gamma_i\|_F^2=\sum_i\|\Lambda_i-\Lambda_iS^{-1}\|^2_F+\sum_i\|\Lambda_iS^{-1}-\Gamma_i\|_F^2.$$
\end{thm}
\begin{pf} We have \begin{align*}
\sum_i\|\Lambda_i-\Gamma_i\|_F^2&=\sum_i\sum_{k=1}^n\|\Lambda_ie_k-\Gamma_ie_k\|^2\\
&=\sum_{k=1}^n\sum_i\|\Lambda_ie_k-\Gamma_ie_k\|^2
\end{align*}
From Proposition \ref{prop1}, we get
\begin{align*}
\sum_{k=1}^n\sum_i\|\Lambda_ie_k-\Gamma_ie_k\|^2&=\sum_{k=1}^n\sum_i\|\Lambda_ie_k-\Lambda_iS^{-1}e_k\|^2
                                                +\sum_{k=1}^n\sum_i\|\Lambda_iS^{-1}e_k-\Gamma_ie_k\|^2\\
                                                &=\sum_i\|\Lambda_i-\Lambda_iS^{-1}\|_F^2+\sum_i\|\Lambda_iS^{-1}-\Gamma_i\|_F^2.
                                                \end{align*}
\end{pf}
\begin{cor}
 Let $\{\Lambda_i\}$ be a $g$-frame for $\mathcal{H}$ finite dimensional Hilbert space, with frame
operator $S$. For all alternate dual g-frames $\{\Gamma_i\}$ of $\{\Lambda_i\}$, the inequality
$$\sum_i\|\Lambda_i-\Gamma_i\|^2_F\geq\sum_i\|\Lambda_i-\Lambda_iS^{-1}\|^2_F$$
takes place and we have equality if and only if $\Gamma_i=\Lambda_iS^{-1}$
\end{cor}
\begin{thm}If $\{\Lambda_i\}$ is a nearly Parseval $g$-frame and $\mathcal{H}$ is a finite n-dimensional Hilbert space, then there exists a dual $\{\Gamma_i\}$ of $\{\Lambda_i\}$ such that the following inequality takes place
$$\sum_i\|\Lambda_i-\Gamma_i\|^2_F\leq n\frac{\varepsilon^2}{1-\varepsilon}.$$
Moreover, the estimation is optimal.
\end{thm}
\begin{pf} We have \begin{align*}\|\Lambda_i-\Lambda_iS^{-1}\|_F^2&=\sum_i\sum_k\|\Lambda_ie_k-\Lambda_iS^{-1}e_k\|^2\\
                                               &=\sum_k\sum_i\|\Lambda_ie_k-\frac{1}{\lambda_k}\Lambda_ie_k\|^2\\
                                               &=\sum_k\bigg(1-\frac{1}{\lambda_k}\bigg)^2\sum_i\|\Lambda_ie_k\|^2\\
                                               &=\sum_k\bigg(1-\frac{1}{\lambda_k}\bigg)^2\langle Se_k,e_k\rangle\\
                                               &=\sum_k\bigg(1-\frac{1}{\lambda_k}\bigg)^2\lambda_k\\
                                               &=\sum_k\bigg(\sqrt{\lambda_k}-\frac{1}{\sqrt{\lambda_k}}\bigg)^2=\sum_k\frac{(\lambda_k-1)^2}{\lambda_k}.
\end{align*}
But, as before, $1-\varepsilon\leq\lambda_k\leq 1+\varepsilon.$ It follows $|\lambda_k-1|\leq\varepsilon$ and $\displaystyle \frac{1}{\lambda_k}\leq\frac{1}{1-\varepsilon}$ and thus $$\sum_i\|\Lambda_i-\Lambda_iS^{-1}\|_F^2\leq n\frac{\varepsilon^2}{1-\varepsilon}.$$
\end{pf}

\textbf{Remarks.} 1. Some results of this paper were presented at the $24^{th}$ International Conference on Operator Theory, July 2-7, 2012, West University of Timi\c soara, Romania.\\
 2. Some results of this paper are related to the ones presented in the papers \cite{Eldar} and \cite{Frank}. The papers \cite{Eldar} and \cite{Frank} deals only with vector frames. Our results are more general and we use a technique which is more simple even in the case of vector frames. In addition, our paper contains results concerning the best approximation of operator-valued frames by its alternate duals, with optimal estimates.\\

\textbf{Acknowledgements.} The authors would like to thank the referees and the editor of this paper for their valuable comments and remarks.

The final work of P. G\u avru\c ta on this paper  was supported by a grant of Romanian National Authority for Scientific Research, CNCS-UEFISCDI, project number PN-II-ID-JRP-2011-2/11-RO-FR/01.03.2013.

\end{document}